\documentclass[12pt]{article}

 \usepackage{amsmath}
 \usepackage{amsfonts}
 \usepackage{amssymb}
 \usepackage{latexsym}
\newtheorem{theorem}{Theorem}

\newtheorem{remark}[theorem]{Remark}
\newtheorem{definition}[theorem]{Definition}

 \allowdisplaybreaks

\begin{document}

\begin{center}

{\Large \bf
A Compactness Theorem for
Homogenization of Parabolic Partial Differential Equations}\\

\bigskip

{\sc By~~Hee Chul Pak } \\
\medskip
{\it Department of Mathematics, Yonsei University }
\end{center}

\vskip 1cm

\begin{abstract}
{\small
In order to have a better description of homogenization for
parabolic partial differential equations with periodic coefficients,
we define the notion of parametric two-scale convergence.
A compactness theorem is proved
to justify this notion.}
\end{abstract}


%%=============================================================================
\section{Introduction}
%%=============================================================================
It is well known that the modeling of physical processes in
strongly inhomogeneous media leads to the study of differential
equations with rapidly varying coefficients.
Regarding  coefficients as periodic functions,
many attempts for getting approximate solutions are accomplished,
and some of successful ones are G-convergence by Spagnolo,
H-convergence by Tartar  and
$\Gamma$-convergence by De Giorgi. Another common way is to use
formal asymptotic expansion  -
we first guess by a formal expansion what the limit should be
and then justify it by energy method.
Two-scale convergence defined by G. Allaire (\cite{1})
is an efficient way of combining these two procedures -
 two-scale convergence guarantees the least amount of
convergent degree, which is stronger than weak convergence and
weaker than norm convergence. But it is restricted to the elliptic
cases. For studying parabolic differential equations, the
convergent nature is not clear and even two-scale convergence can
not be applied directly.
Even though stationary problems corresponding to the
given parabolic problems may be considered for homogenization process,
the convergent relationship
between stationary problems and non-stationary
problems is not   justified
by the two-scale convergence defined by G. Allaire.
In this note,
    we present a new approach which is a modification of
    the two-scale convergence  technique
to explain  the homogenization of
    parabolic partial differential equations
    with periodic coefficients, and this also describes the convergence
nature of parabolic equations on inhomogeneous media.
This notion is justified by a compactness theorem which is the main
result of this paper.

%%==================================================================
\section{Parametric two-scale convergence}
%%==================================================================

Let $\Omega$ be an open set in ${\large\bf R}^n$ and
$Y \equiv [0,1]^n$.
We denote various spaces of periodic functions by a subscript $\sharp$.
    For example,
         $C_{\sharp}^{\infty}(Y)$ denotes the space of infinitely
        differentiable functions on ${\bf R}^n$
        that are periodic
        of period $Y$.

We define the parametric two-scale convergence for parabolic differential
equations with periodic coefficients following the lead of
Nguetseng \cite{4}.
\begin{definition}
A sequence of functions $u^{\varepsilon}(t, x) \in
 L^{\infty}([0,\infty); L^2(\Omega))$
   is said to parametric two-scale converge  to a limit
$U_0(t,x,y) \in
                    L^{\infty}([0,\infty);L^2(\Omega \times Y))$
             if for any test function
             $\Psi(t, x,y) \in L^1\left([0,\infty); C_0^{\infty}
           (\Omega ;C_{\sharp}^{\infty}(Y))\right)$,
we have
\[
\lim_{\varepsilon \rightarrow 0}
      \int_0^{\infty}\int_{\Omega}
                          u^{\varepsilon}(t,x)
                          \Psi \! \left(t,x,\frac{x}{\varepsilon}
                   \right)
                      dx
       dt
    = \int_0^{\infty}\int_{Y} \int_{\Omega}
                                   U_0(t,x,y)\Psi(t,x,y)
                              dx
                     dy
       dt.
\]                                                        %                                                             %
\end{definition}          %%&&&&&&&&&&&&&&&&&&&&&&&&&&&&&&&&&&&&&&&&

    Roughly speaking, the convergence is a kind of
    $weak^*$-convergence with respect to
    $L^{\infty}(0,\infty)$-norm and two-scale with respect to
    $L^2(\Omega)$-norm.

\begin{remark}
In \cite{1}, we can find the following;
 For any $\Psi(x,y) \in L^2(\Omega;C_{\sharp}(Y))$,
    we have
\[
\lim_{\varepsilon \rightarrow 0} \int_{\Omega} \Psi \!
\left(x,\frac{x}{\varepsilon} \right)^2dx = \int_{\Omega}\int_{Y}
\Psi(x,y)^2dxdy.
\]
\end{remark}

We prove a compactness theorem for the
parametric two-scale convergence.

\begin{theorem}
Any bounded sequence $ u^{\varepsilon}$
 in $L^{\infty}([0,\infty); L^2(\Omega))$ has a
 parametric two-scale convergent subsequence.
\end{theorem}
\noindent {\bf Proof : }  %%****************************
Suppose that  $ u^{\varepsilon}$
     is a bounded sequence in
     $L^{\infty}([0,\infty); L^2(\Omega))$.
We want to show that there is a subsequence
           $ u^{\varepsilon_j}$ of
           $ u^{\varepsilon}$ and a function
           $U_0(t,x,y) \in
                    L^{\infty}([0,\infty);L^2(\Omega \times Y))$
such that
for any $\Psi(t,x,y)$  in
          $L^1\left([0,\infty); C_0^{\infty}
           (\Omega ;C_{\sharp}^{\infty}(Y))\right) $
\[
\lim_{\varepsilon_j \rightarrow 0}
      \int_0^{\infty}\int_{\Omega}
                          u^{\varepsilon_j}(t,x)
                          \Psi \! \left(t,x,\frac{x}{\varepsilon_j}
                   \right)
                      dx
       dt
    = \int_0^{\infty}\int_{Y} \int_{\Omega}
                                   U_0(t,x,y)\Psi(t,x,y)
                              dx
                     dy
       dt.
\]
Let ${\cal F}_{\varepsilon}(\Psi)
         \equiv
               \int_0^{\infty}\int_{\Omega}
                      u^{\varepsilon}(t,x)
                      \Psi(t,x,\frac{x}{\varepsilon})
                dx dt $
 and
        ${\cal D} \equiv L^2(\Omega; C_{\sharp}(Y))$.
There is a positive constant $C$ such that
                                 $ \| u^{\varepsilon}
                                   \|_{L^{\infty}(0, \infty;
                                                L^2(\Omega))} \leq C$.
Since
\begin{eqnarray*}
\mid {\cal F}_{\varepsilon}(\Psi) \mid
    &\leq & \int_0^{\infty}  \left\|  u^{\varepsilon}(t)
                             \right\|_{L^2(\Omega)}
                    \left\|\Psi \! \left(t,x,\frac{x}{\varepsilon} \right)
                    \right\|_{L^2(\Omega)}
             dt \\
    &\leq & C \int_0^{\infty} \max_{y \in Y}
            \parallel\Psi(t,x,y)\parallel_{L^2(\Omega)} dt
            =    C \int_0^{\infty}\parallel\Psi(t)\parallel_{\cal D} dt
\end{eqnarray*}
           and
        $\Psi \in L^1\left([0,\infty) ; {\cal D} \right) $.
We have ${\cal F}_{\varepsilon} \in
                             \left(L^1([0,\infty) ; {\cal D})\right)'$.
    So a subsequence ${\cal F}_{\varepsilon_j}$
       is $weak^{*}$-convergent
       to some ${\cal U}_0
                 \in
                 \left(L^1(0,\infty ; {\cal D})\right)'$.
Hence for any
       $\Psi$ in $L^1([0,\infty) ; {\cal D})$,
 we have
\[
\left\|\Psi \! \left(t,x,\frac{x}{\varepsilon_j} \right)
\right\|
         _{L^2(\Omega)}
\; \leq \;
           \parallel \Psi(t)
           \parallel_{{\cal D} },
\]
so by Lebesgue Dominated Convergence Theorem,
\begin{eqnarray*}
\mid{\cal U}_0(\Psi)\mid
             =   \lim_{\varepsilon_j \rightarrow 0}
               \mid {\cal F}_{\varepsilon_j}(\Psi) \mid
    &\leq &  \limsup_{\varepsilon_j \rightarrow 0}
             \int_0^{\infty} \parallel
                         u^{\varepsilon_j}(t)
             \parallel_{L^2(\Omega)}
             \left\|\Psi \! \left(t,x,\frac{x}{\varepsilon_j} \right)
             \right\|
                      _{L^2(\Omega)} dt \\
    &\leq &  C \int_0^{\infty}
               \limsup_{\varepsilon_j \rightarrow 0}
               \left\| \Psi \! \left(t,x,\frac{x}{\varepsilon_j} \right)
               \right\|
                      _{L^2(\Omega)} dt   \\
     &  = &   C \int_0^{\infty}
               \parallel \Psi(t)
               \parallel_{ L^2(\Omega \times Y)}dt.
\end{eqnarray*}
Since $L^1([0,\infty) ; {\cal D}) $ is dense in $L^1([0,\infty);
       L^2(\Omega \times Y))$,
      it follows  that  ${\cal U}_0$ is in
      $\left(L^1([0,\infty) ; L^2(\Omega \times Y))
       \right)' $.
By Riesz Representation Theorem,
    we have
\[
   {\cal U}_0(\Psi)
             = \int_0^{\infty} \left< U_0(t), \Psi(t) \right>
               _{L^2(\Omega \times Y)} dt
\]
for some $U_0 \in L^{\infty}
          \left([0,\infty); L^2(\Omega \times Y)
          \right)$.
Therefore
\begin{eqnarray*}
\lim_{\varepsilon_j \rightarrow 0}
    \int_0^{\infty} \!\!\int_{\Omega} \!
              u^{\varepsilon_j}(t,x)
              \Psi \! \left(t,x,\frac{x}{\varepsilon_j} \right)
              dx dt
  &\equiv&    \lim_{\varepsilon_j \rightarrow 0}
              {\cal F}_{\varepsilon_j}\!(\Psi)        \\
  &  =   &    {\cal U}_0(\Psi)                         \\
  &  =   &    \int_0^{\infty}\!\!\!
              \int_{Y} \!\int_{\Omega} \!
              U_0(t,x,y)\Psi(t,x,y)dxdydt.
\end{eqnarray*}
\hfill$\Box$\par
  %%****************************************************************

\section*{Acknowledgements}
This paper was supported
by Korea Research Foundation, Grant No. KRF-2001-005-D20004.

%%%%%%%%%%%%%%%%%%%%%%%%%%%%%%%%%%%%%%%%%%%%%%%%%%%%%%%%%%%%%%%%%%%%%%%%%%%%%%
\noindent Department of  Mathematics, Yonsei University, \\
134 Sinchon-dong, Seodaemun-gu, Seoul 120-749, Korea
\\E-mail: hpak@yonsei.ac.kr


\bigskip


\begin{thebibliography}{99}
{\footnotesize
%
\bibitem[1]{1} %%{allaire}
 G. Allaire,
Homogenization  and two-scale convergence, SIAM J. Math. Anal.
23:1482 - 1518, 1992.
%
\bibitem[2]{2} %%{lion}
 A.Benoussan, J. L. Lions and G.
Papanicolaou, Asymptotic Analysis for Periodic Structures.
North Holland, Amsterdam, 1978.
%
\bibitem[3]{3} %%{hornung}
U. Hornung,
Homogenization and Porous Media,
Springer, Berlin-New York, 1996.
%
\bibitem[4]{4}
G. Nguetseng, A general convergence result for a functional
related to the theory of homogenization, SIAM J. Math. Anal.,
20 (1989), pp. 608 - 623.
%
\bibitem[5]{5} %%{heechul} %%{show}
H-C Pak, R.E. Showalter,  Thin-Film Capacitance
Models, Applicable Analysis, To appear.
%
\bibitem[6]{6} %%{show}
R. E. Showalter,
 Monotone Operators in Banach Space
and Nonlinear Partial Differential Equations,
volume 49 of
Mathematical Surveys and Monographs, AMS, 1997.
%
\bibitem[7]{7} %%{zko}
V.V. Zhikov, S.M. Kozlov, and O.A. Oleinik,
Homogenization of Differential Operators and Integral Functionals,
Springer-Verlag, Berlin-New York, 1991.
%
}


\end{thebibliography}
\end{document}